\documentclass[12pt,a4paper]{article}

\usepackage[dvipdfmx]{color}

\usepackage[pdftex]{graphicx}

\usepackage{authblk}

\usepackage{amsmath}
\usepackage{amsthm}
\usepackage{ascmac}
\usepackage{amssymb}
\usepackage{amsfonts}
\usepackage{mathtools}
\usepackage{bm}
\usepackage{here}
\usepackage{cancel}
\usepackage{comment}
\usepackage{fancyhdr}
\usepackage{cases}
\usepackage{multicol}
\usepackage{braket}
\usepackage{geometry}
\theoremstyle{definition}

\newtheorem*{theorem*}{Theorem}

\newtheorem*{definition*}{Definition}

\newtheorem*{remark*}{Remark}

\newtheorem*{proposition*}{Proposition}

\newtheorem*{lemma*}{Lemma}

\title{Shogi and Frieze group}
\author{Yusuke Imai}
\affil{Graduate School of Information Science and Technology, The University of Tokyo, 7-3-1 Hongo, Bunkyo-ku, Tokyo 113-8656, Japan}
\date{}

\begin{document}
\maketitle

\begin{abstract}
Shogi is a traditional Japanese strategy board game in the same family as chess, chaturanga, and xiangqi, and has been theoretically studied from various aspects. 
The research on recommended sequences of moves in each opening of shogi is called \emph{joseki}; how to use a rook (Static Rook and Ranging Rook), or how to develop a castle, etc.
Also, many pieces of \emph{tsume shogi}, artistic shogi miniature problems, in which the opponent's king is checkmated by a series of checks, have been created involving various beautiful techniques such as ``saw" and ``puzzle ring".
In addition, the rapid development of AI in recent years has led to the pursuit of the best possible moves in shogi.
In this paper, we move away from the study of winning and losing in shogi and focus on the mathematical aspects of the movement of shogi pieces.
We propose to correspond movements of shogi pieces to a set of geometrical patterns constructed by the shape of shogi pieces and representing the Frieze group through the condition regarding the neighborhood of arrangements of given shogi pieces.
Although the discovery of this correspondence does not lead to a winning strategy for shogi, however, it does demonstrate a curious involvement between the traditional Japanese board game and Western mathematics.
\end{abstract}
 


\section{Introduction}
Shogi is a traditional Japanese strategy board game  that has much in common with chess, chaturanga, and xiangqi.
According to the website of Japan Shogi Association (\emph{Nihon Shogi Renmei}) \cite{JSA}, the most plausible explanation for the origin of Shogi is that it originated in the ancient Indian game of chaturanga, and it is not known when shogi was born.
The \emph{Nichu Reki}, an encyclopedia of customs estimated to have been compiled between 1210 and 1221, includes a prototype of shogi called Heian shogi [Heian small shogi (\emph{Heian sho shogi}) using 6 kinds of pieces and Heian large shogi (\emph{Heian dai shogi}) using 13 kinds of pieces], which shows the influence of Chaturanga pieces.
Around the 13th century, large shogi (\emph{dai shogi}) that is an extension of the Heian shogi with more pieces began to be played, and small shogi (\emph{sho shogi}) was invented, which incorporated the \emph{hisha}, \emph{kakugyo}, and \emph{suizo} of \emph{dai shogi} into Heian shogi. 
In the 15th and 16th century, the \emph{suizo} was removed from \emph{sho shogi}, and the game became the present \emph{hon shogi}.

Theoretical studies of shogi can be divided into two main categories: \emph{joseki} and \emph{tsume shogi}.
The \emph{joseki} is the study of the openings of shogi while \emph{tsume shogi} is miniature problems related to the closings of shogi.
There are two aspects to \emph{joseki}, one concerning the position of rook (\emph{hisha}) and the other concerning the castle of king (\emph{gyokusho}).
The strategy of not moving rook sideways at the beginning of the game is called Static Rook (\emph{ibisha}), while the strategy of moving rook sideways is called Ranging Rook (\emph{furibisha}).
The Static Rook includes many strategies: \emph{kaku gawari, ai gakari, yoko hu dori}, and so on. 
The Ranging Rook also includes many strategies: \emph{shiken bisha, sanken bisha, naka bisha}, and so on.
In addition, there are many types of the castle of king: \emph{mino gakoi, yagura gakoi, anaguma gakoi}, and so on.
Each of them is meticulously constructed and further subdivided, and in some cases, the study continues until the checkmate.
The \emph{tsume shogi} has various additional rules not found in regular shogi and requires a high level of mathematical ability in its creation, however it is basically a game of checkmating the opponent's king with a series of checks.
The \emph{tsume shogi} was already created in the Edo period, and a collection of \emph{tsume shogi} pieces by the first master (\emph{shodai meijin}), Ohashi Sokei I, titled \emph{shogi zobutsu} was written in 1602 \cite{Tsume}.
There are various types of \emph{tsume shogi} depending on types of the techniques, and some of them were already shown in the collections of artistic \emph{tsume shogi} pieces from the Edo period \emph{Shogi Muso} by the seventh master (\emph{nanasei meijin}) Ito Sokan III and \emph{Shogi Zuko} by Ito Kanju, including ``saw of dragon horse move" (\emph{uma noko}), ``smoke checkmate" (\emph{kemuri zume}), and ``dragon chase" (\emph{ryu oi}) with 611 checks (In Refs. \cite{Tsume1,Tsume2,Tsume3}, you can perform the procedure of them).
Kurushima Yoshihiro, a researcher of mathematics (\emph{wasan}) in the Edo period, made early references to Euler's totient function and cofactor expansion \cite{Kurushima}, while also developing ``puzzle ring" in \emph{tsume shogi} \cite{Tsume4}.
Since then, various \emph{tsume shogi} techniques have been developed, and many \emph{tsume shogi} pieces are still published in ``Tsume Shogi Paradise" \cite{Tsumepara}.
In recent years, with the development of AI, research on shogi incorporating AI has been active. 
The prominent example is that the Shogi AI named ``elmo" won the Masuda Kozo Prize in 2020, which is awarded to those who have made excellent new moves or contributed to the advancement of \emph{joseki} \cite{elmo}.

The study of shogi so far has involved a systematic study of the openings and a study of the closings with high artistry and sophisticated mathematical techniques.
This paper sheds light on a new aspect of shogi research by focusing solely on the movement of each shogi piece, away from winning and losing.
We show that a correspondence exists between the shogi piece's movement and the geometric patterns made by shape of shogi pieces and characterized by Frieze groups. 
 
\newpage
\section{Movement of shogi pieces}
\begin{figure}[H]
\centering
{%
\resizebox*{\textwidth}{!}{\includegraphics{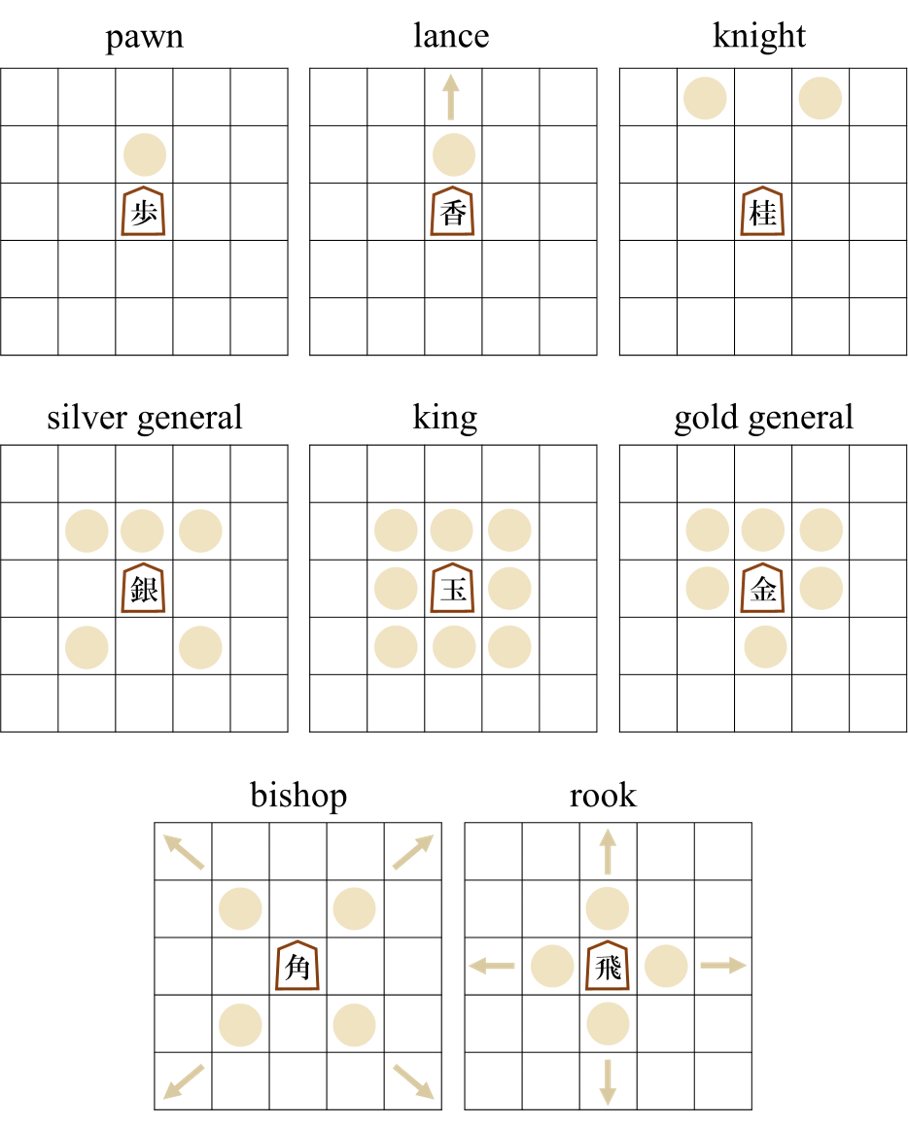}}}\hspace{5pt}
\caption{Movement of the shogi pieces: the pawn (\emph{fuhyo}), lance (\emph{kyosha}), knight (\emph{keima}), silver general (\emph{kinsho}), king (\emph{osho/gyokusho}), gold general (\emph{kinsho}), bishop (\emph{kakugyo}), and rook (\emph{hisha}) from left-top to right-bottom. The disks show the possible movements of each shogi piece and the arrow means that the shogi piece can advance in that direction as long as it does not collide with another piece.} \label{Fig1}
\end{figure}

\begin{figure}[H]
\centering
{%
\resizebox*{\textwidth}{!}{\includegraphics{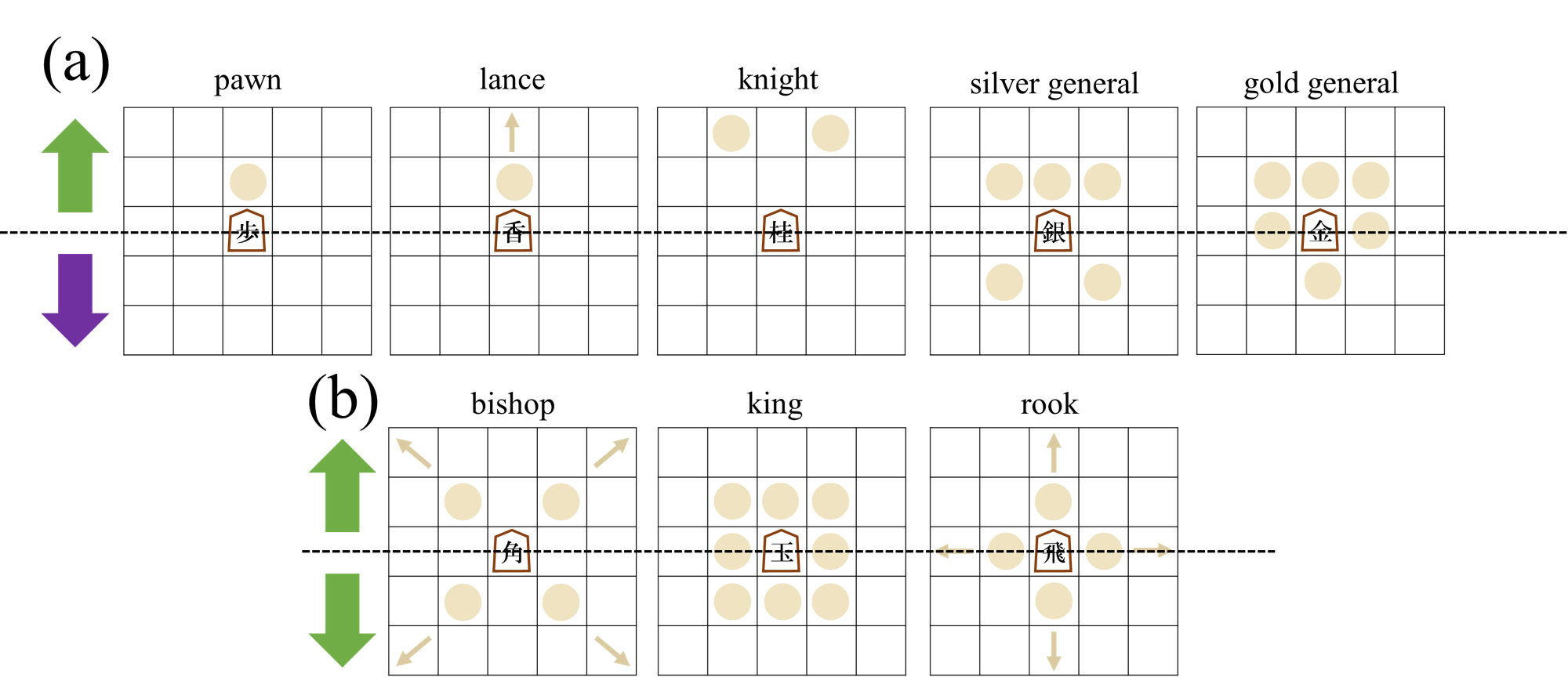}}}\hspace{5pt}
\caption{(a) The movements of the shogi pieces with the horizontal mirror symmetry and (b) without the horizontal mirror symmetry.} \label{Fig2}
\end{figure}

We first remark that the movements of shogi pieces are complicated and do not contain simple mathematical structures.
One can see the following two points from Fig. 1 showing the movement of all the shogi pieces.
\\
\\
1. Many shogi pieces have moves that break horizontal mirror symmetry.
\\
\noindent
2. Not all pieces allowed by symmetry within a certain range are included in the set of shogi pieces.
\\
\\
The first statement is obvious from Fig. \ref{Fig1}; the movement of the king, rook, and bishop have horizontal symmetry [Fig. \ref{Fig2}(a)] while the movement of the gold general, silver general, knight, lance, and pawn do not have horizontal symmetry [Fig. \ref{Fig2}(b)].
Note that the chess pieces have horizontal mirror symmetry except for pawns.
The second statement, for example, concerns the movement of the lance.
While the lance can move forward as many squares as possible, it cannot move backward at all.
A piece that moves backward any number of squares but not forward at all is not included in the set of shogi pieces.
Note that the set of shogi pieces does not include pieces that have symmetrical up-and-down movement of the lance, but such a piece does exist in the big shogi (\emph{dai shogi}), and it is called reverse chariot (\emph{hensha}).
Thus, the movements of shogi pieces are complicated, and their mathematical structure cannot be simply comprehended.

\newpage
\section{Shogi pattern}
\subsection{Shogi pattern and neighborhood}
\begin{figure}[H]
\centering
{%
\resizebox*{0.25\textwidth}{!}{\includegraphics{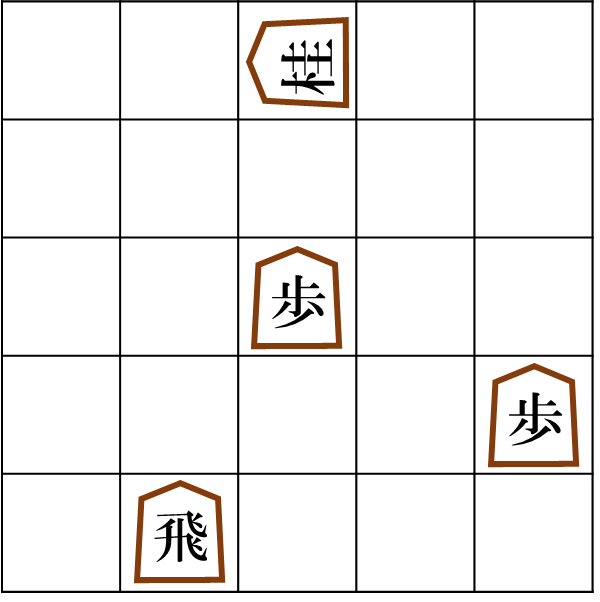}}}\hspace{5pt}
\caption{The example of a shogi pattern.} \label{Fig3}
\end{figure}

The situation discussed in the previous section would lead us away from trying to understand the mathematical structure of a shogi piece's movement by the movement of a single shogi piece alone.
In this paper we consider geometric patterns constructed by some kinds of shogi pieces, shogi patterns, as shown in Fig. \ref{Fig3}.
Shogi patterns are composed of what ${\bf form}$ the shogi pieces make and by what ${\bf kinds}$ ${\bf of}$ ${\bf shogi}$ ${\bf pieces}$ they are made.
Note that shogi pieces are not placed horizontally in a regular game, but in this paper they are assumed to be placed horizontally, as is done in \emph{yonin shogi} (shogi for four players).
We assume that pieces of the same orientation cannot be taken from each other.
This paper mainly deals with only the simple shopi patterns constructed by one kind of shogi pieces.

\begin{figure}[H]
\centering
{%
\resizebox*{\textwidth}{!}{\includegraphics{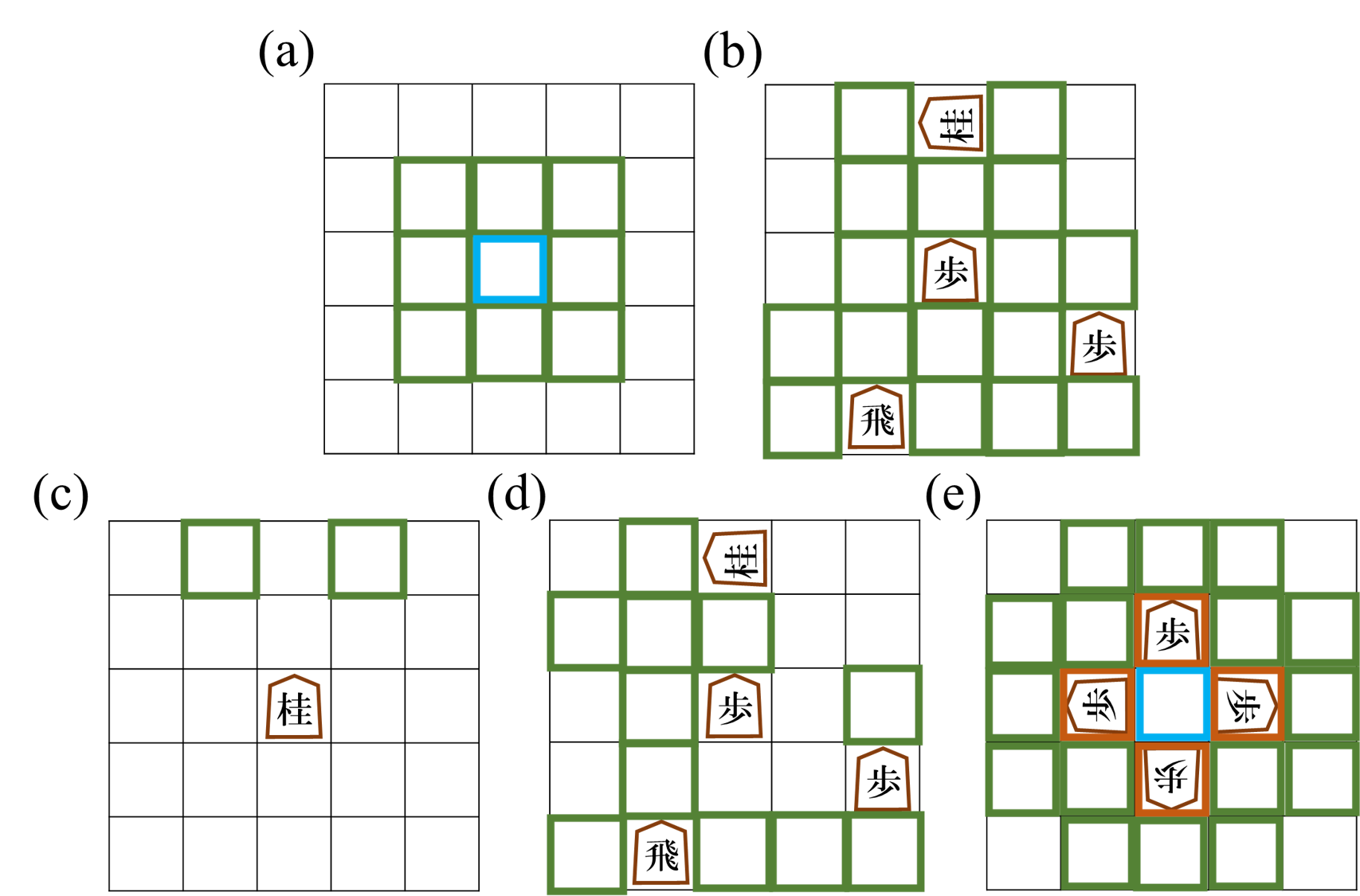}}}\hspace{5pt}
\caption{The example of (a) neighborhood (green) of a given square (blue), (b) neighborhood of a shogi pattern, (c) the control of a shogi pieces, (d) the control of a shogi pattern, (e) inside (blue), base (brown), outside (green) neighborhood of a shogi pattern.} \label{Fig4}
\end{figure}

Here we make some definitions concerning the shogi pattern.
First, we define the neighborhood of a given square is a set of squares that can be placed vertically, horizontally, and diagonally in a single square from the given square [Fig. \ref{Fig4}(a)]. 
Then, we define a neighborhood of a given shogi pattern is all the squares included in all the neighborhoods of the squares included in the given shogi pattern [Fig. \ref{Fig4}(b)].
We define the control of a given shogi piece as the squares where the piece can move [Fig. \ref{Fig4}(c)].
We define the control of a given shogi pattern as the squares where the pieces included in the given shogi pattern can move [Fig. \ref{Fig4}(d)].
We define an inside neighborhood of a given shogi pattern as squares that is included in the neighborhood of the shogi pattern and surrounded by shogi pieces [Fig. \ref{Fig4}(e)].
In addition we define a base neighborhood of a given shogi pattern as the set of squares constructing the shogi pattern and included in the neighborhood of the given shogi pattern [Fig. \ref{Fig4}(e)].
We also define an outside neighborhood of a given shogi pattern as the set of squares that is included in the neighborhood of the shogi pattern and not included in the inside neighborhood and in the base neighborhood [Fig. \ref{Fig4}(e)].

\subsection{Neighoborhood control condition}
Then, we introduce the following condition for a square included in the neighborhood of a given shogi pattern.
\\
\\
\noindent
[${\bf Neighoborhood}$ ${\bf condrol}$ ${\bf condition}$] 
\\
\noindent
For a square included in the neighborhood of a given shogi pattern, the square satisfies the neighborhood control condition if and only if the square is included in the control of the given shogi pattern. 
\\
\\
Here we define some properties for a shogi pattern.
A shogi pattern satisfies ${\bf complete}$ ${\bf neighborhood}$ ${\bf control}$ ${\bf condition}$ if and only if all the squares included in neighborhoods of the shogi pattern satisfy the neighborhood control condition.
Shogi patterns consisting of shogi pieces with one direction cannot satisfy the above condition if there exists the base neighborhood in the shogi patterns since the allied pieces cannot be taken, thus weaker conditions must be considered.
Then we introduce the following condition; a shogi pattern satisfies ${\bf nearly}$ ${\bf compelte}$ ${\bf neighborhood}$ ${\bf control}$ ${\bf condition}$ if and only if a part of the neighborhood of the shogi pattern is included in the control and the set of squares that are included in the neighborhood of the shogi pattern and that do not satisfy the neighborhood control condition is either the empty set (i.e., the complete neighborhood condition is satisfied), base neighborhood, inside neighborhood, or outside neighborhood.

\begin{figure}[H]
\centering
{%
\resizebox*{\textwidth}{!}{\includegraphics{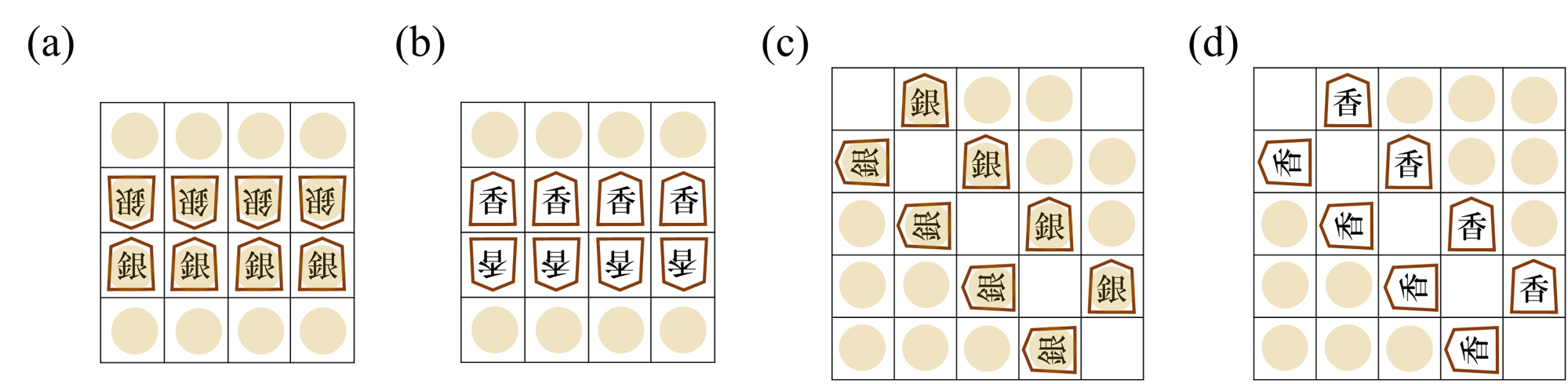}}}\hspace{5pt}
\caption{The example of a shogi pattern that (a) satisfies the complete neighborhood control condition, (b), (c) satisfies the nearly complete neighborhood control condition, (d) does not satisfy the complete neighborhood control condition or nearly complete neighborhood control condition.} \label{Fig5}
\end{figure}

For example, Fig. \ref{Fig5}(a) shows a shogi pattern that satisfies the complete neighborhood control condition, and Figs. \ref{Fig5}(b) and (c) show a shogi pattern that satisfies the nearly complete neighborhood control condition and Fig. \ref{Fig5}(d) shows a shogi pattern that does not satisfy the complete neighborhood control condition and nearly complete neighborhood control condition.
Although the actual shogi board is 9 squares by 9 squares, we assume that the shogi board and the shogi patterns in Fig. \ref{Fig5} and in the following figures are infinitely extended.

In the following subsections, we study some properties relating to the complete neighboring control condition and nearly complete neighboring control condition.

\subsection{Specialty of king}
We remark that any shogi pattern constructed by the kings placed so that they are not next to each other always satisfies nearly complete neighborhood control condition because the king can move to all the neighboring squares (see Fig. \ref{Fig1}).
This property is unique to the king and no other kinds of shogi pieces have this property.

\subsection{Special shogi form}
\begin{figure}[H]
\centering
{%
\resizebox*{\textwidth}{!}{\includegraphics{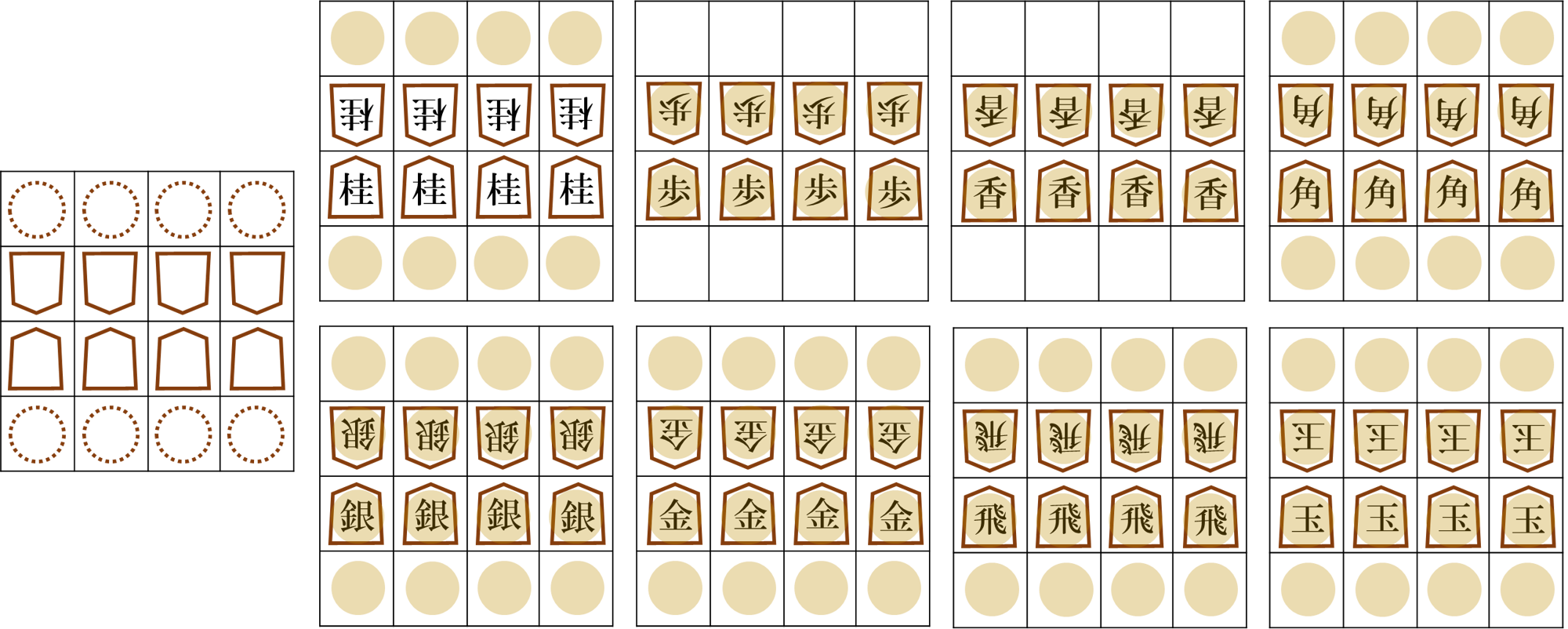}}}\hspace{5pt}
\caption{An example of the special shogi forms where the dotted circles represent the outside neighborhood.} \label{Fig6}
\end{figure}

We remark that there exist some special shogi forms such that for any kind of shogi pieces, the shogi pattern constructed by them satisfies the nearly complete neighborhood control condition.
For example, a shogi form shown in the left of Fig. \ref{Fig6} is such a shogi form.
The outside (base) neighborhood of the shogi pattern constructed by the knight is (not) included in the control.
The base (outside) neighborhood of the shogi pattern constructed by the pawn or lance is (not) included in the control.
The base and outside neighborhood of the shogi pattern constructed by the other kinds of shogi pieces are included in the control.

\subsection{Dualities}
\begin{figure}[H]
\centering
{%
\resizebox*{0.75\textwidth}{!}{\includegraphics{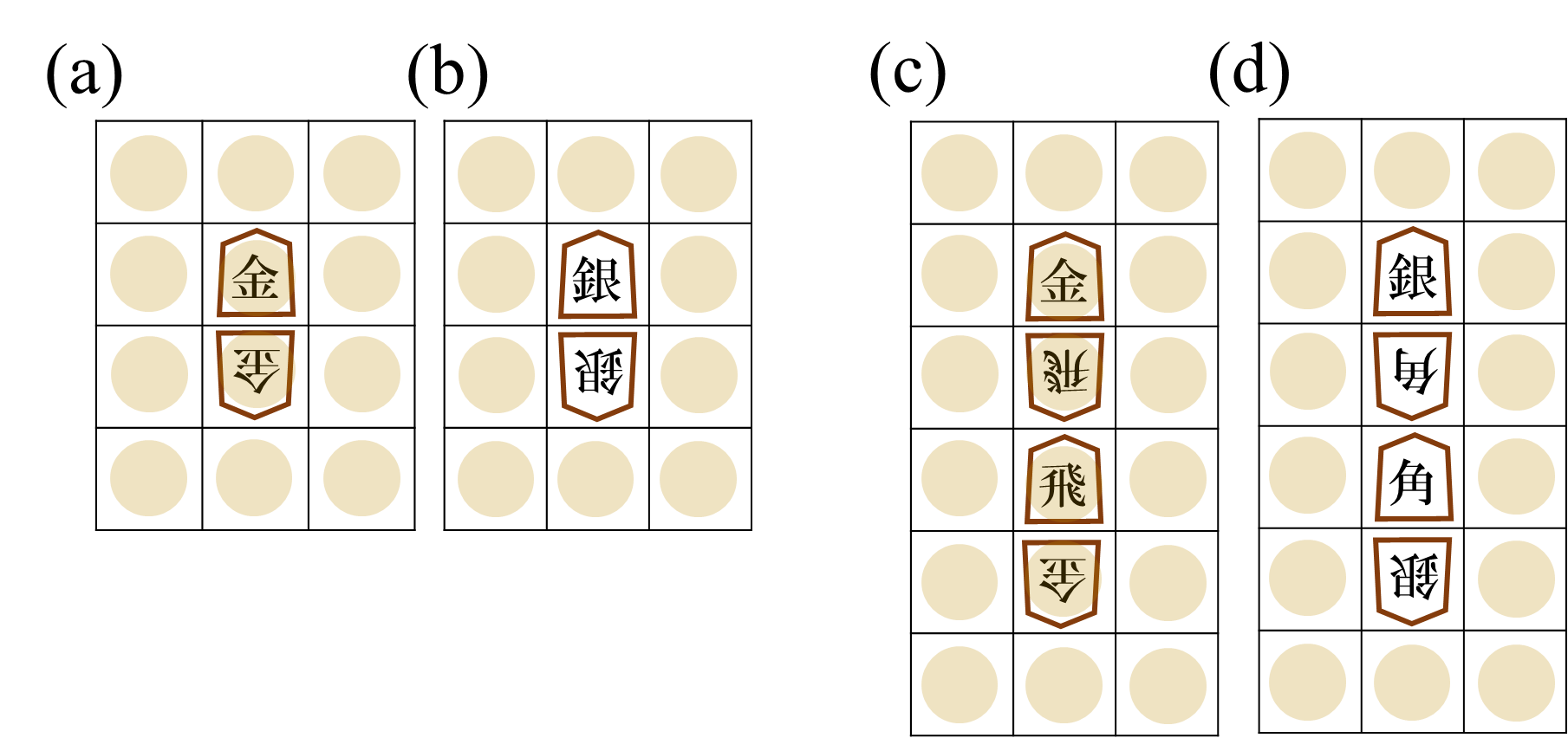}}}\hspace{5pt}
\caption{(a), (b) A duality between the gold general and silver general through the nearly complete neighborhood control condition and complete neighborhood control condition. (c), (d) A duality between the gold-general/rook and silver-general/bishop through the nearly complete neighborhood control condition and complete neighborhood control condition.} \label{Fig7}
\end{figure}

Figures \ref{Fig7}(a) and (b) show shogi patterns that show a duality between the gold general and silver general through the nearly complete neighborhood control condition and complete neighborhood control condition.
Figure \ref{Fig7}(a) is constructed by the gold general and satisfies the complete neighborhood control condition and Fig. \ref{Fig7}(b) is constructed by the silver general and satisfies the nearly complete neighborhood control condition.
In addition, Figs. \ref{Fig7}(c) and (d) show shogi patterns that show a duality between the gold-general/rook and silver-general/bishop through the nearly complete neighborhood control condition and complete neighborhood control condition.
Figure \ref{Fig7}(c) is constructed by the gold general and rook, and satisfies the complete neighborhood control condition and Fig. \ref{Fig7}(d) is constructed by the silver general and bishop, and satisfies the nearly complete neighborhood control condition.

\newpage
\section{Frieze group and shogi crystal}
In this section we review the Frieze group and introduce shogi crystals each of them has each of the symmetries included in the Frieze groups.
The Frieze groups are mathematical concepts representing two-dimensional patterns that repeat in one direction \cite{Frieze}.
The Frieze groups are a set of 7 groups constructed by translations, glide reflections, vertical reflections, horizontal reflections, and $180^\circ$ rotations: p1, p11g, p1m1, p11m, p2, p2mm, p2mg.
The p1 contains only translations. 
The p11g contains glide-reflections and translations.
The p1m1 contains vertical reflection lines and translations.
The p11m contains translations, horizontal reflections, and glide reflections.
The p2 contains translations and $180^\circ$ rotations.
The p2mg contains vertical reflections, glide reflections, translations and $180^\circ$ rotations.
The p2mm contains horizontal and vertical reflections, translations and $180^\circ$ rotations.

\subsection{Typical geometrical patterns}
\begin{figure}[H]
\centering
{%
\resizebox*{\textwidth}{!}{\includegraphics{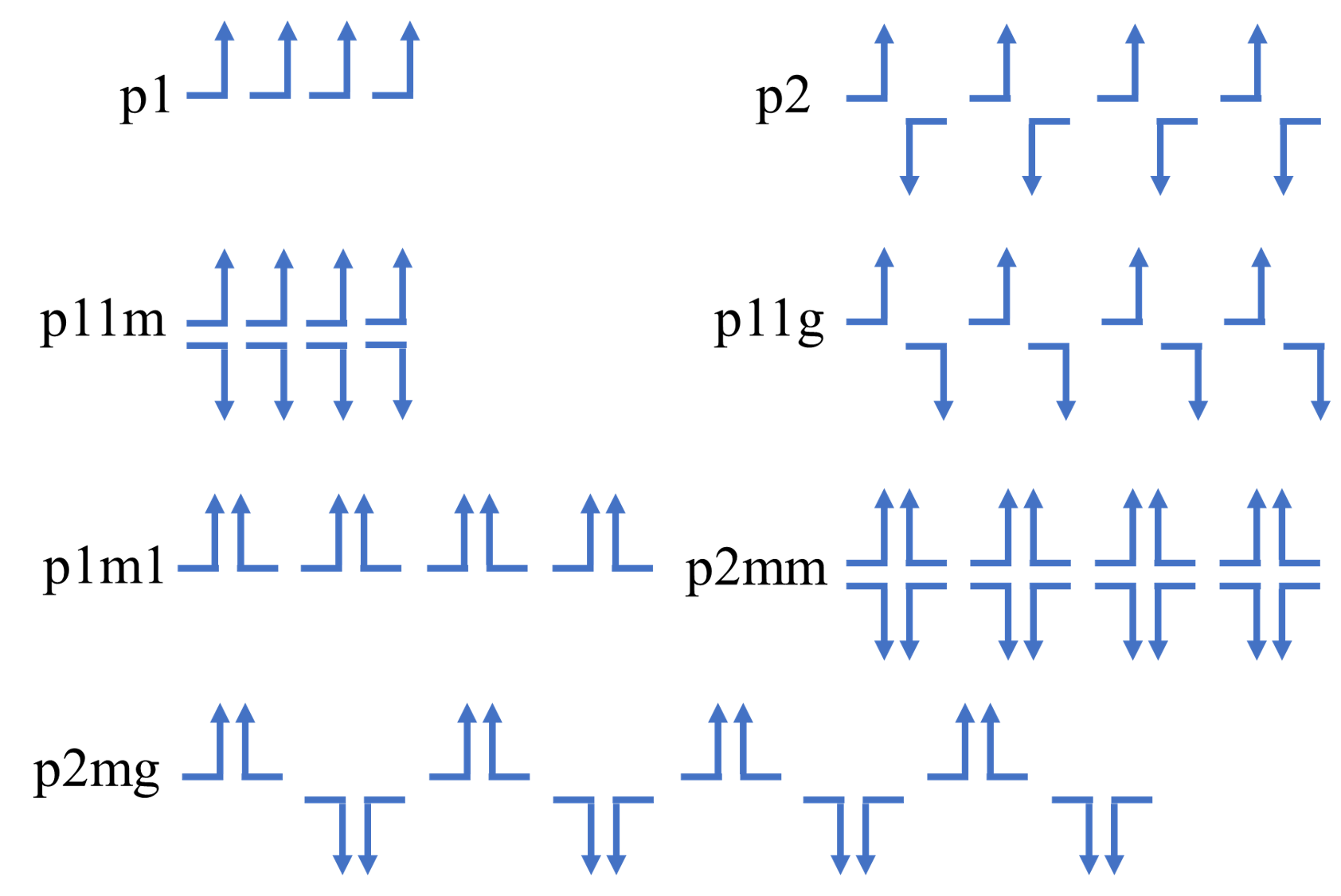}}}\hspace{5pt}
\caption{A simple example of geometrical patterns representing the Frieze groups.} \label{Fig9}
\end{figure}

\begin{figure}[H]
\centering
{%
\resizebox*{\textwidth}{!}{\includegraphics{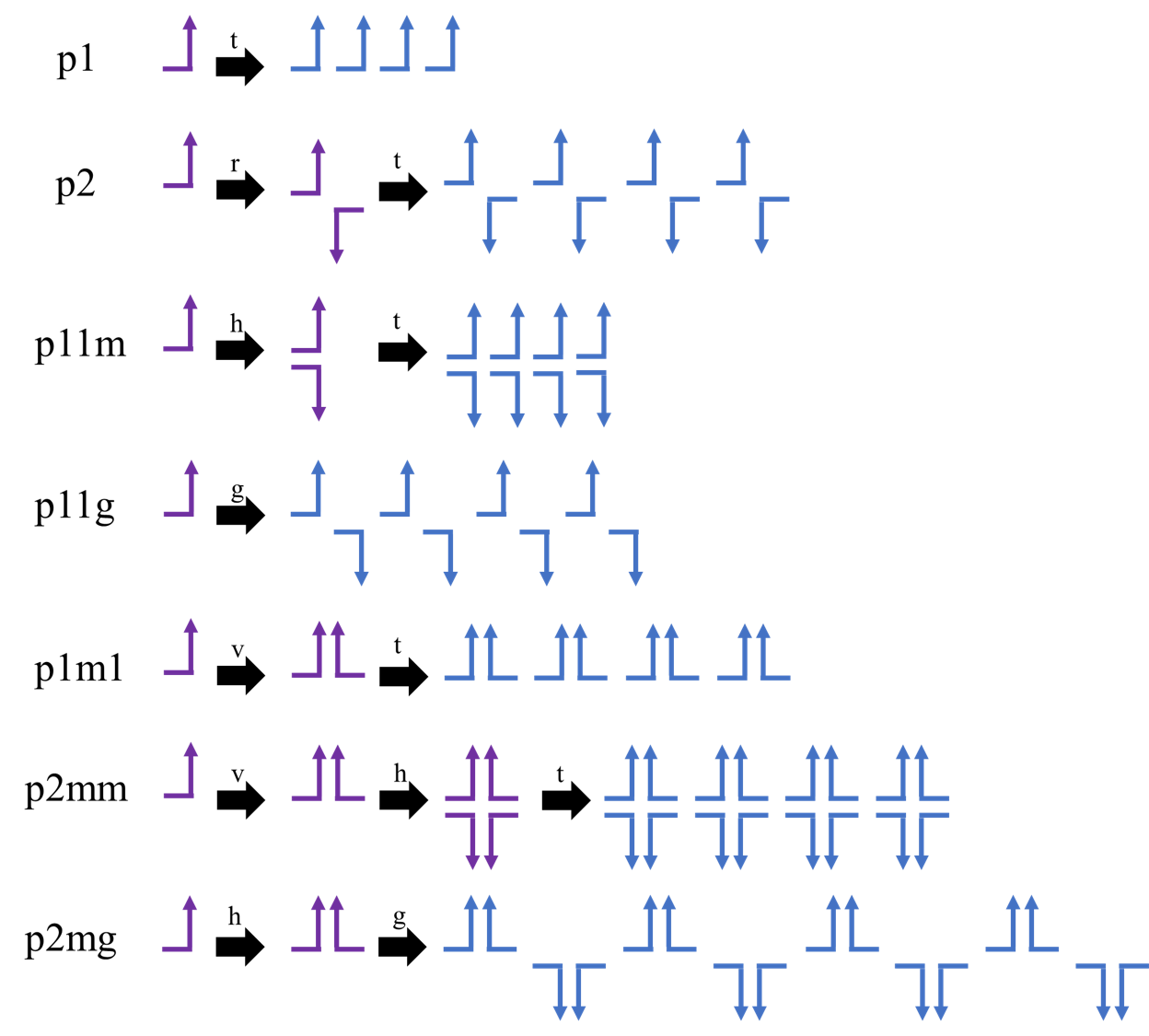}}}\hspace{5pt}
\caption{A figure showing how each of the figures in Fig. \ref{Fig7} is constructed where ``t" represents translation, ``r" represents $180^\circ$ rotation, ``g" represents glide reflection, ``h" represents horizontal reflection, and ``v" represents vertical reflection.} \label{Fig10}
\end{figure}

Simple geometrical patterns representing the Frieze groups are summarized in Fig. \ref{Fig9}.
Here, we have considered a figure of ``L" inverted left and right with an arrow above it as the basic figure for the figures shown in Fig. \ref{Fig9}.
If another set of figures is considered instead by changing the basic figure to something else without any symmetries, the exact same argument about symmetry holds true.
Figure \ref{Fig10} shows how each of the figures in Fig. \ref{Fig9} is constructed and one can see the symmetry of each figure clearly.

\subsection{Shogi crystal}
\begin{figure}[H]
\centering
{%
\resizebox*{\textwidth}{!}{\includegraphics{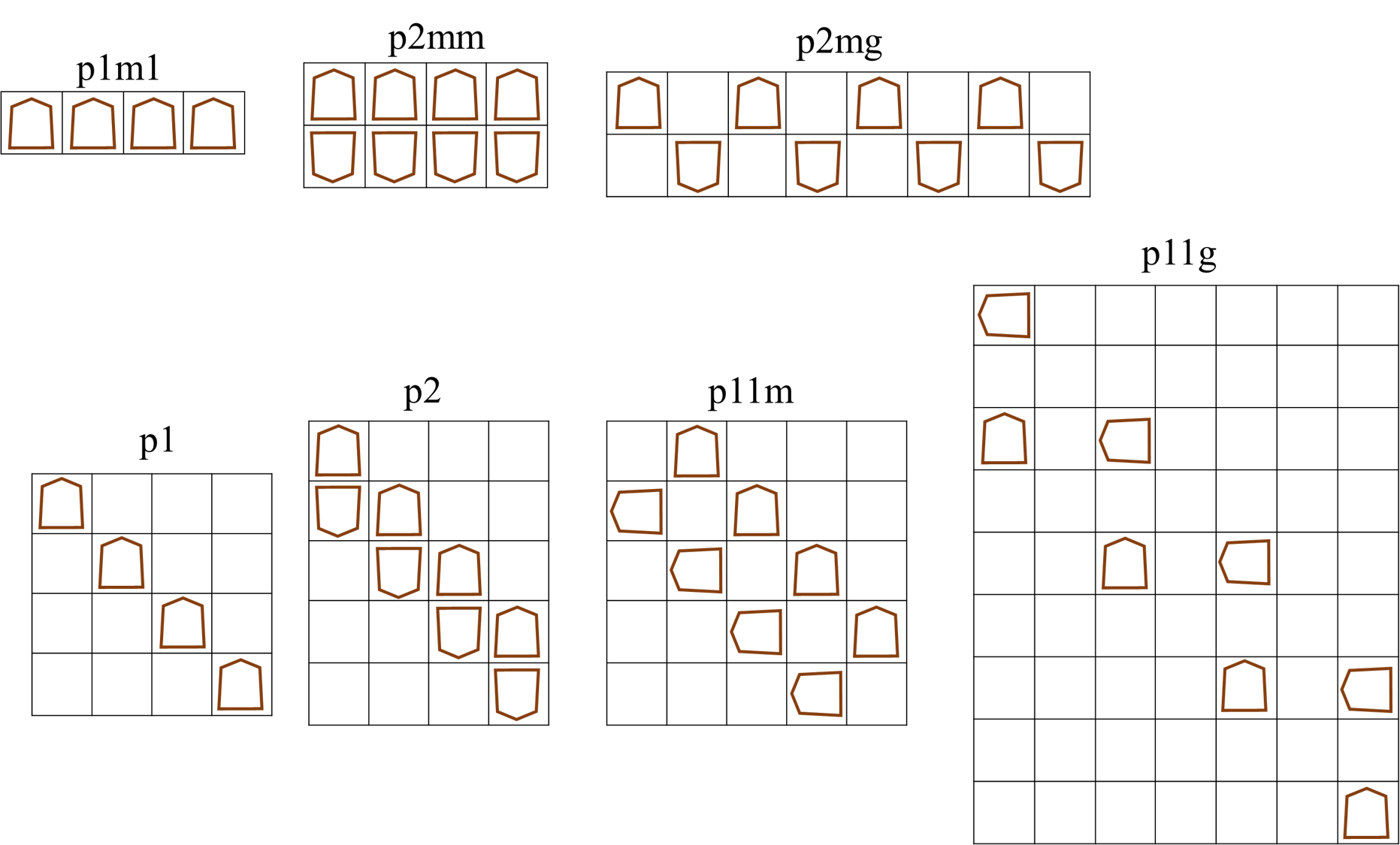}}}\hspace{5pt}
\caption{Shogi crystals.} \label{Fig11}
\end{figure}

\begin{figure}[H]
\centering
{%
\resizebox*{\textwidth}{!}{\includegraphics{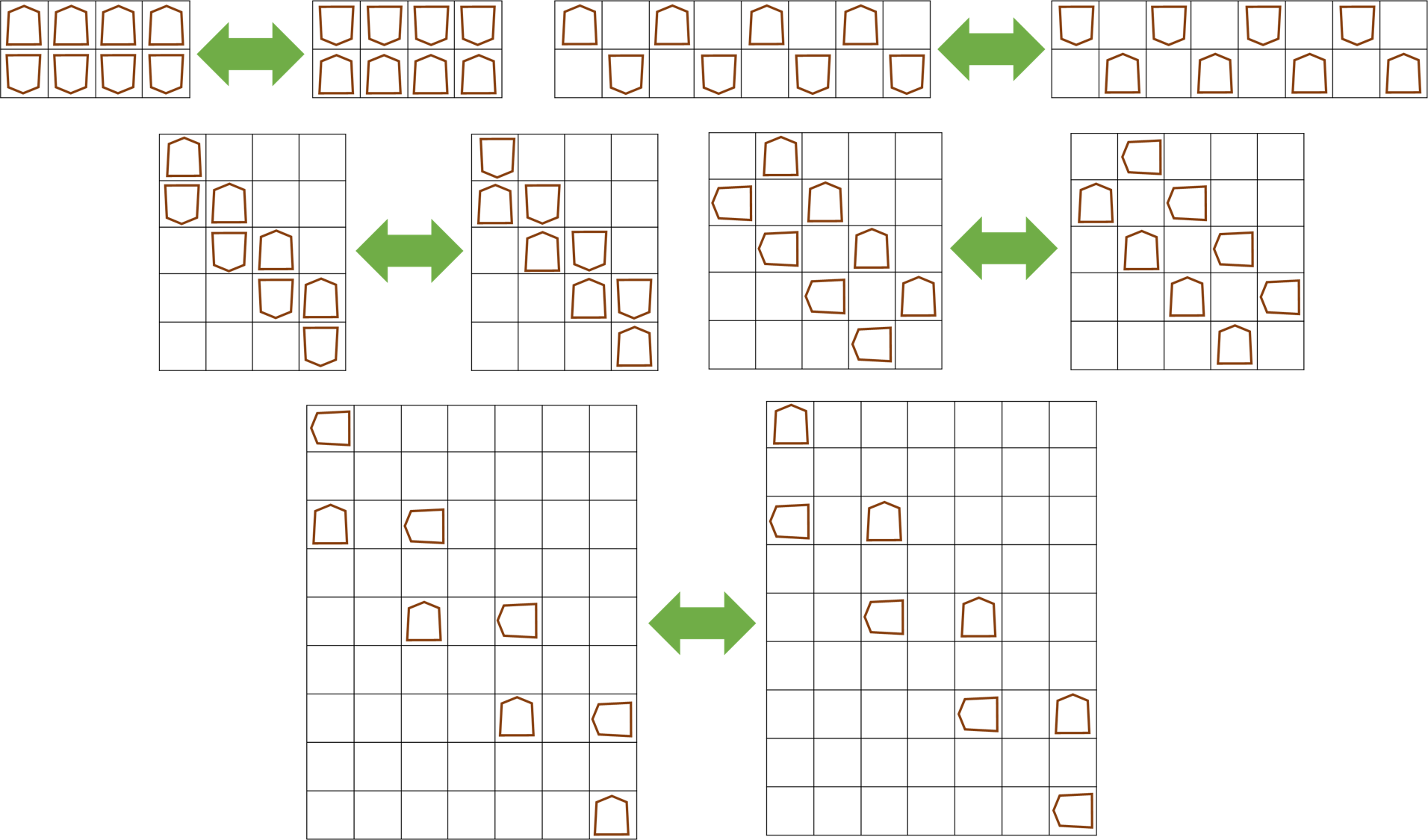}}}\hspace{5pt}
\caption{Shogi crystals and their dual.} \label{Fig12}
\end{figure}

\begin{figure}[H]
\centering
{%
\resizebox*{0.25\textwidth}{!}{\includegraphics{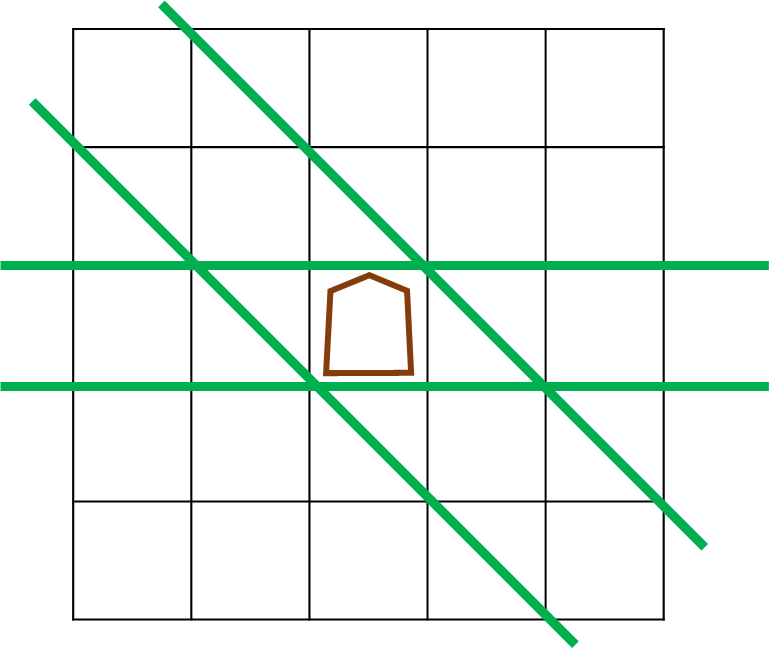}}}\hspace{5pt}
\caption{The mirror axes used in the construction of the shogi crystals.} \label{Fig13}
\end{figure}

In this subsection, we introduce shogi crystals that are similar to figures shown in Fig. \ref{Fig8} using the shape of the shogi pieces.
We remark that the shogi pieces are vertically-symmetrical in shape.
Therefore, for the figures introduced in the previous subsection that are vertically-symmetrical, similar figures can be created using the shape of shogi pieces.
However, for the figures that are introduced in the previous subsection and that are vertically-asymmetrical, it is not possible to create directly similar figures by using the shape of shogi pieces. 
We therefore achieve the vertical asymmetry by changing the direction of travel diagonally.
We show all the shogi crystals in Fig. \ref{Fig11}.
Some shogi crystals consist of only one shogi-piece orientation and the others consist of two shogi-piece orientations. 
As for the latter, the two orientations can be interchanged to create the dual shogi crystals with the same symmetry (Fig. \ref{Fig12}).
With such shogi crystals, the correspondence obtained in the next section can be obtained by swapping the movement of shogi pieces up and down.
In the construction of these shogi crystals, we use only the mirror axes shown in Fig. \ref{Fig13}.
Note that if other mirror axes were used, other shogi crystals could be obtained, and even with the mirror axes used in this paper, it is possible to consider shogi crystals with more complicated forms.

\newpage
\section{Shogi crystal and nearly complete neighborhood control condition}
\begin{figure}[H]
\centering
{%
\resizebox*{\textwidth}{!}{\includegraphics{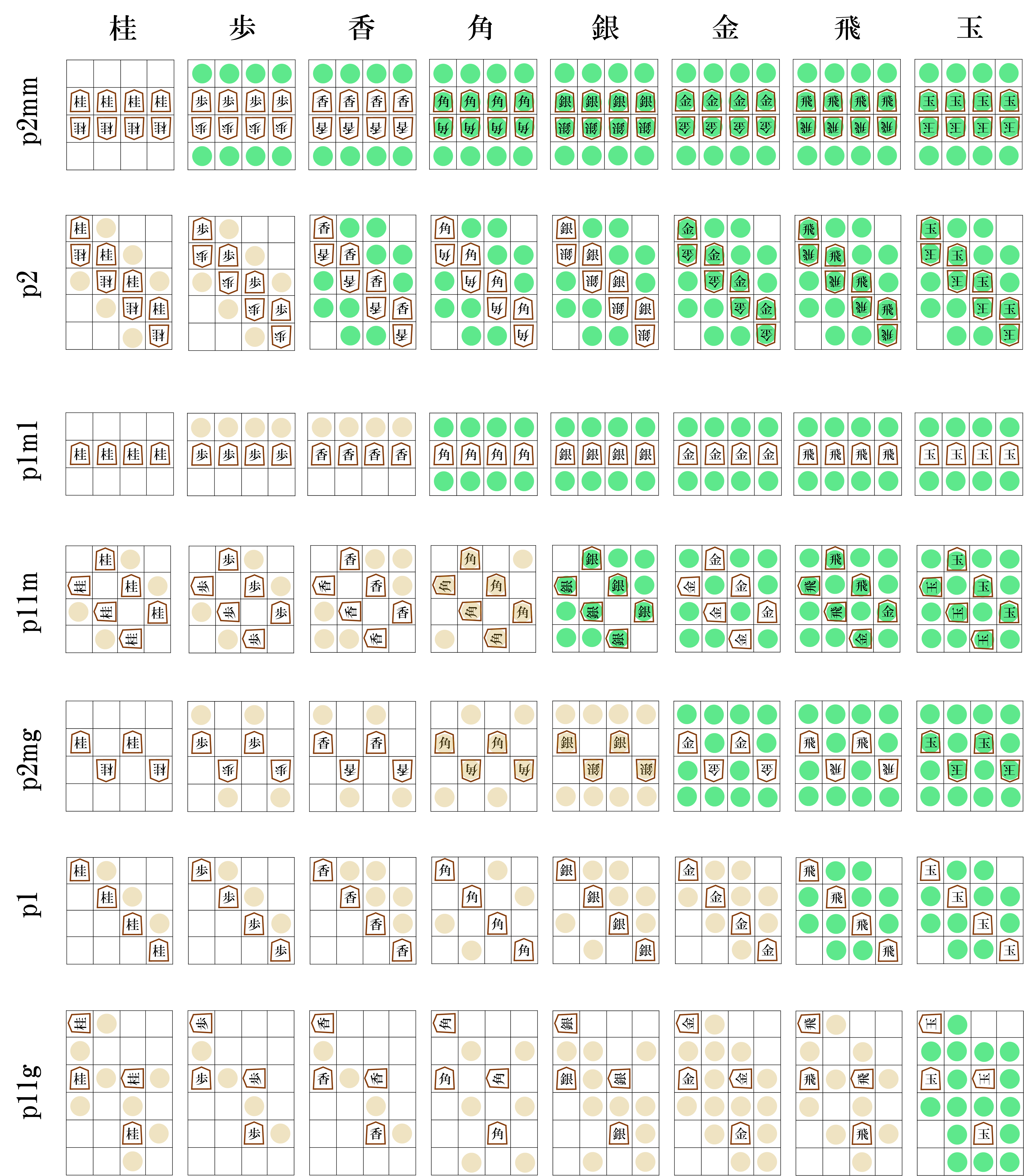}}}\hspace{5pt}
\caption{Summary of neighborhood control of shogi patterns for all the shogi crystals and all the kinds of shogi pieces where the disks show the control of a given shogi pattern and the color of the disks is green if the nearly complete control condition is satisfied.} \label{Fig14}
\end{figure}

\begin{figure}[H]
\centering
{%
\resizebox*{\textwidth}{!}{\includegraphics{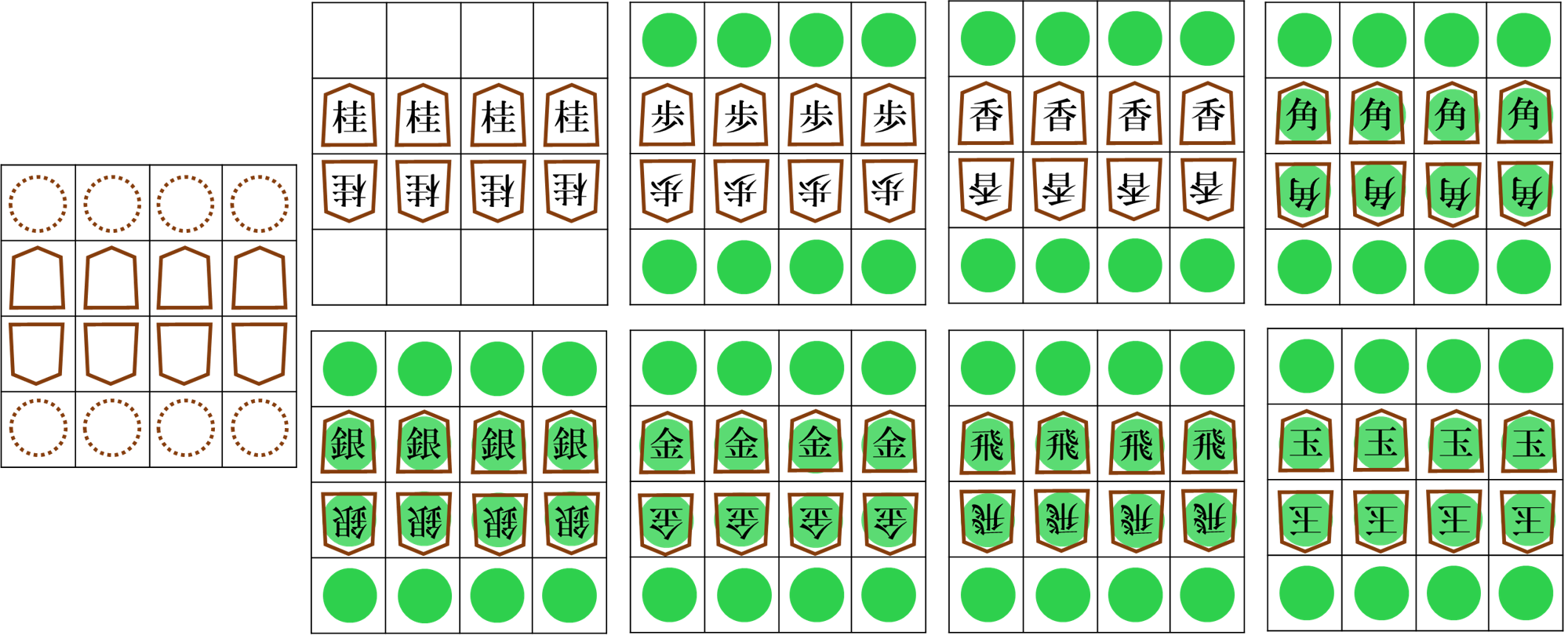}}}\hspace{5pt}
\caption{Neighborhood control of shogi patterns based on the shogi crystal of p2mm symmetry where the green disks show the control in the case that the nearly complete neighborhood control condition is satisfied.} \label{Fig15}
\end{figure}

\begin{figure}[H]
\centering
{%
\resizebox*{\textwidth}{!}{\includegraphics{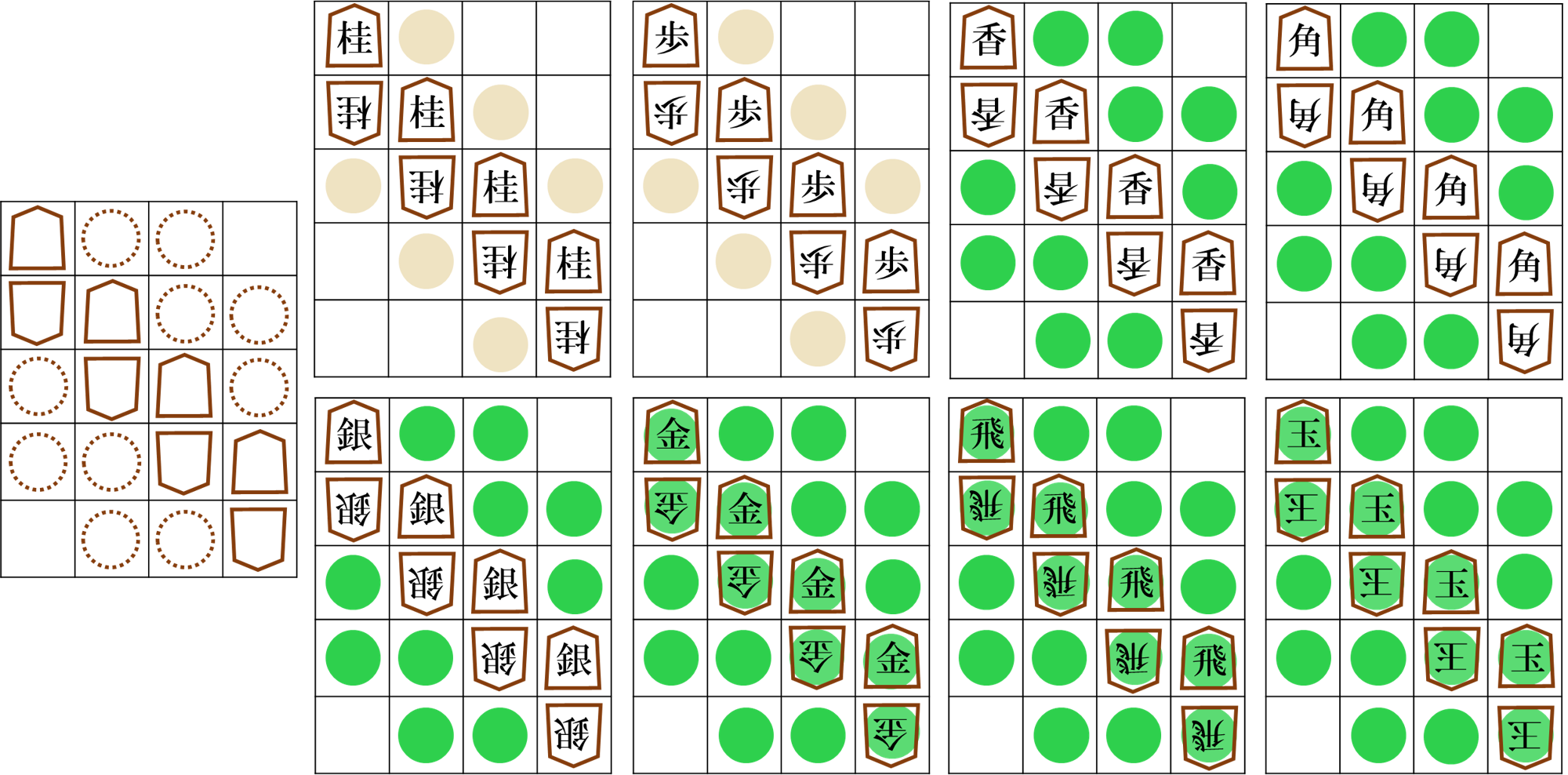}}}\hspace{5pt}
\caption{Neighborhood control of shogi patterns based on the shogi crystal of p2 symmetry where the light gold discs show that a portion of the outside/inside/base neighborhood is included in the control.} \label{Fig16}
\end{figure}

\begin{figure}[H]
\centering
{%
\resizebox*{\textwidth}{!}{\includegraphics{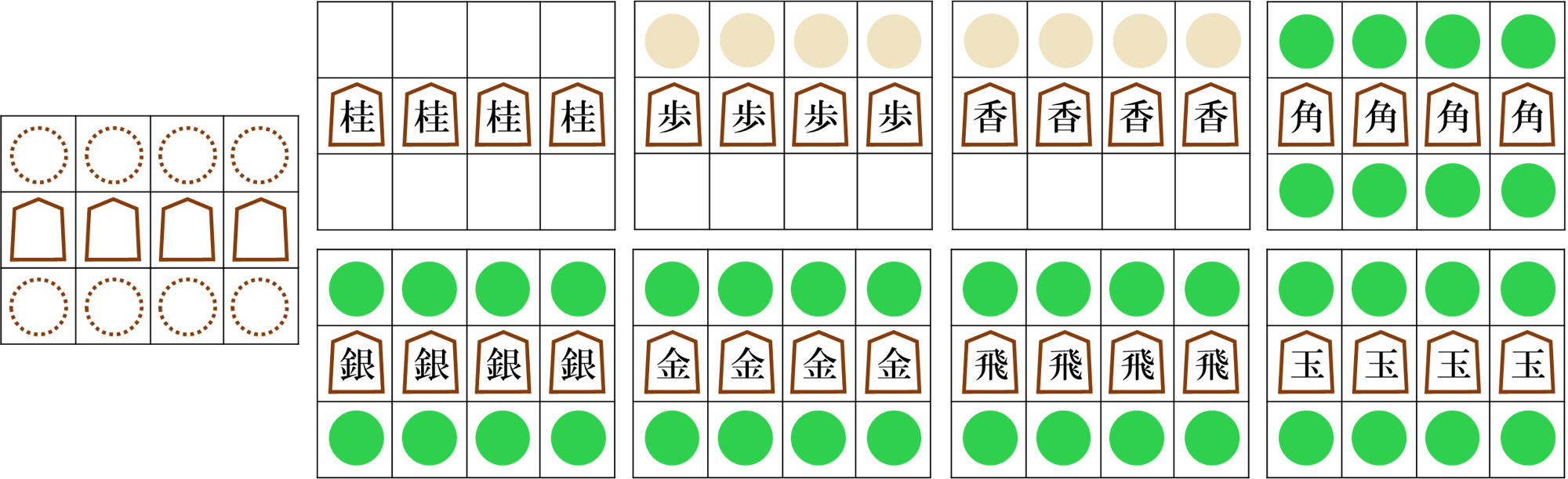}}}\hspace{5pt}
\caption{Neighborhood control of shogi patterns based on the shogi crystal of p1m1 symmetry.} \label{Fig17}
\end{figure}

\begin{figure}[H]
\centering
{%
\resizebox*{\textwidth}{!}{\includegraphics{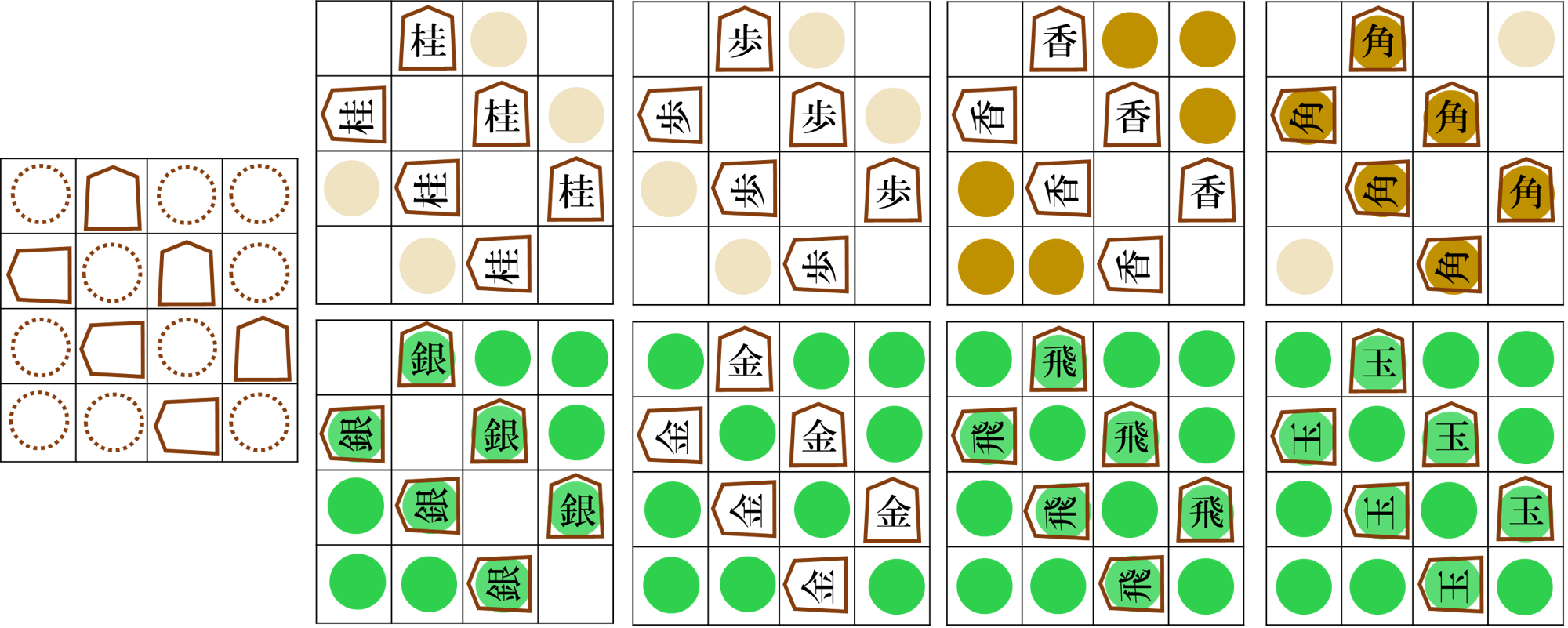}}}\hspace{5pt}
\caption{Neighborhood control of shogi patterns based on the shogi crystal of p11m symmetry where the dark gold discs show that one of the outside/inside/base neighborhood is included in the control.} \label{Fig18}
\end{figure}

\begin{figure}[H]
\centering
{%
\resizebox*{\textwidth}{!}{\includegraphics{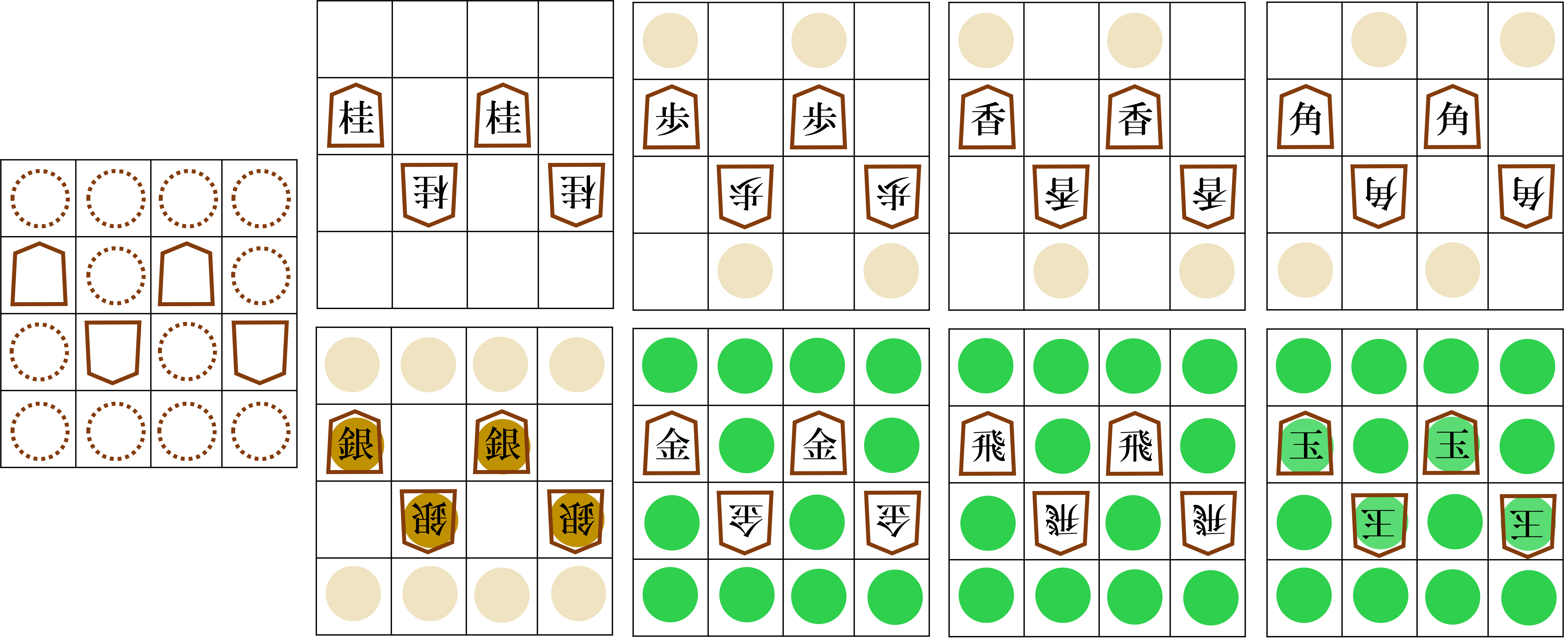}}}\hspace{5pt}
\caption{Neighborhood control of shogi patterns based on the shogi crystal of p2mg symmetry.} \label{Fig19}
\end{figure}

\begin{figure}[H]
\centering
{%
\resizebox*{\textwidth}{!}{\includegraphics{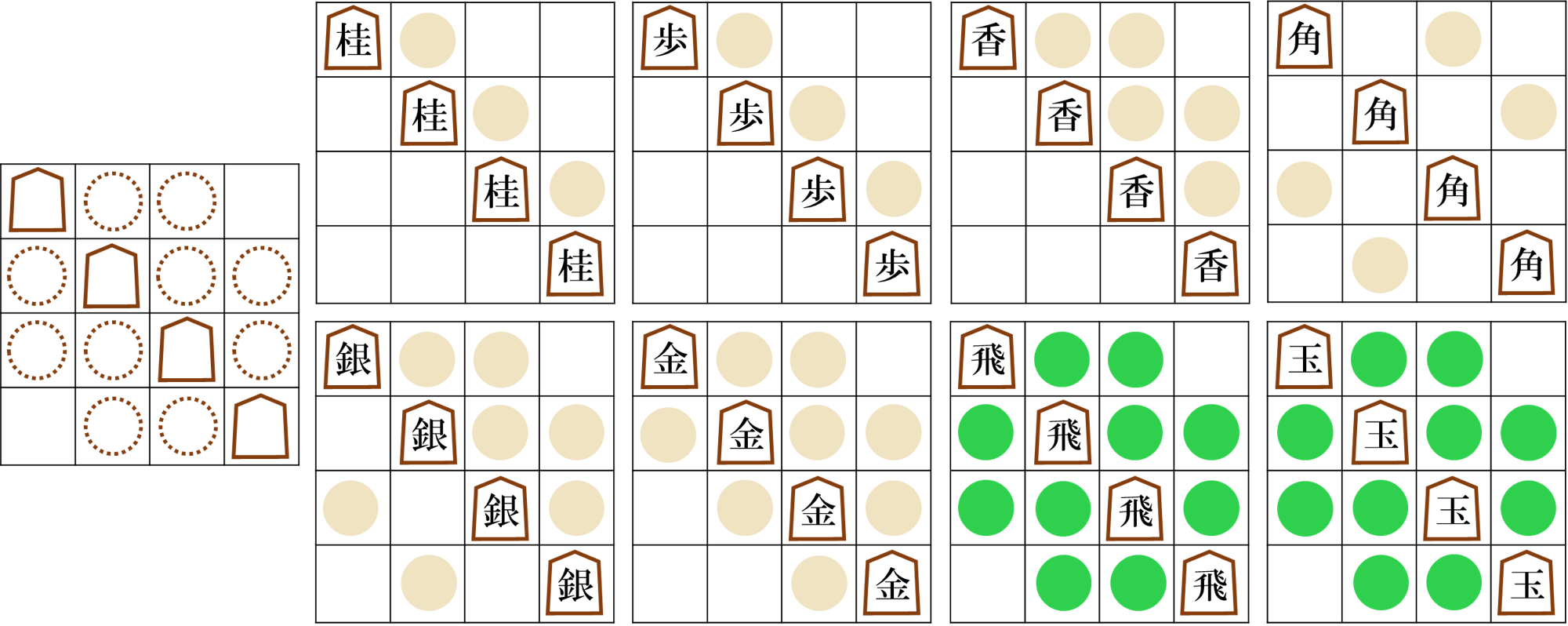}}}\hspace{5pt}
\caption{Neighborhood control of shogi patterns based on the shogi crystal of p1 symmetry.} \label{Fig20}
\end{figure}

\begin{figure}[H]
\centering
{%
\resizebox*{\textwidth}{!}{\includegraphics{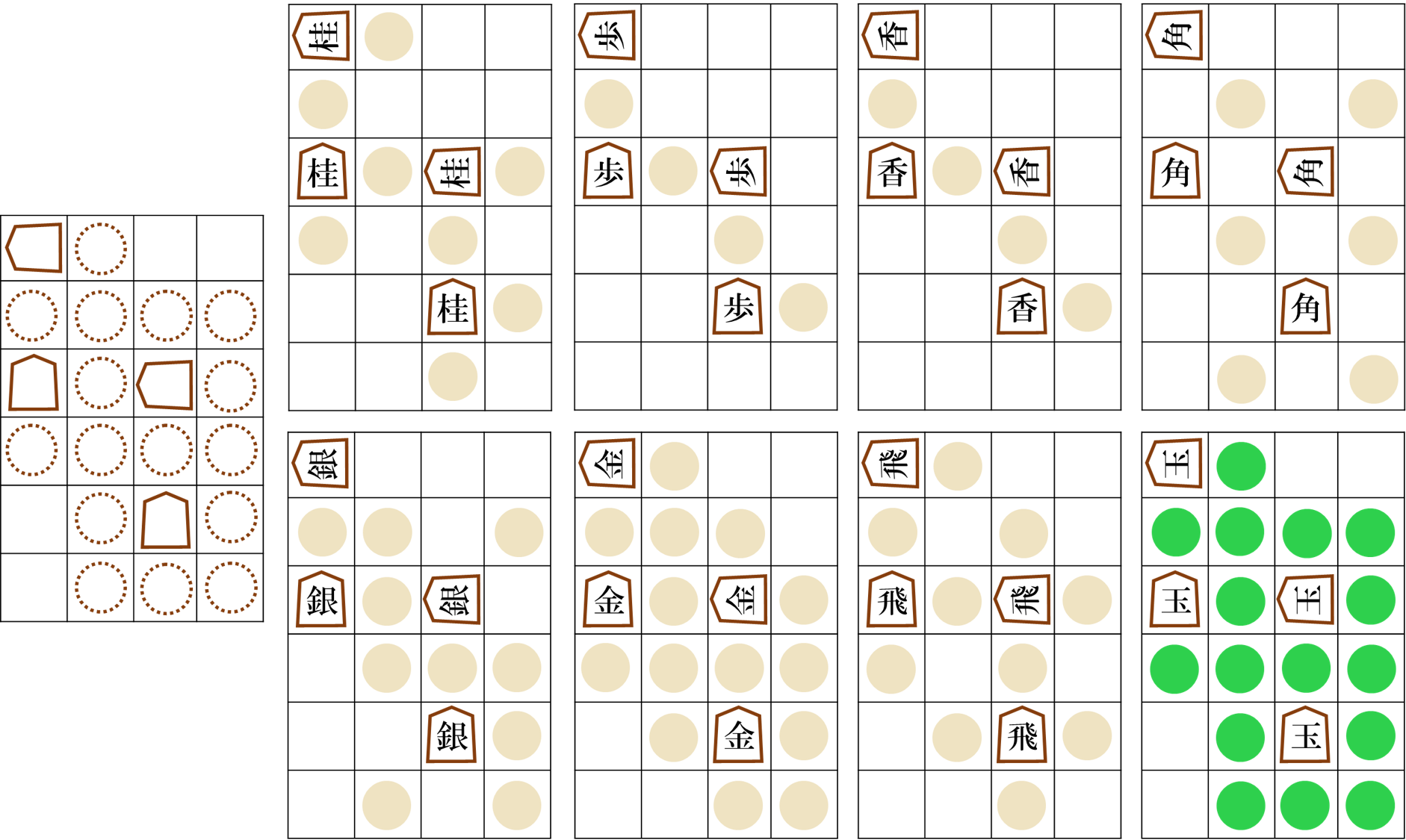}}}\hspace{5pt}
\caption{Neighborhood control of shogi patterns based on the shogi crystal of p11g symmetry.} \label{Fig21}
\end{figure}

In this section, we will examine whether the shogi patterns based on the shogi crystals satisfy the nearly complete neighborhood control condition.
Figure \ref{Fig14} shows the correspondence between the movements of shogi pieces and the Frieze groups through the nearly complete neighborhood control condition.
By using Figs. \ref{Fig15}-\ref{Fig21}, we explain neighborhood control of shogi patterns for all the shogi crystals and all the kinds of shogi pieces in detail.

First we consider the shogi crystal of p2mm symmetry.
As shown in Fig. \ref{Fig15}, the shogi patterns constructed by the shogi crystal of p2mm symmetry and any kinds of shogi pieces except the knight satisfy the nearly complete neighborhood control condition.
Also, the shogi pattern constructed by the shogi crystal of p2mm symmetry and the knight does not satisfy the nearly complete neighborhood control condition.
Note that the knight is the only piece that cannot go forward or forward diagonally.
Also, in the shogi patterns constructed by the shogi crystal of p2mm symmetry and the pawn or lance, the outside (base) neighborhood is (is not) included in their control.
In addition, in the shogi patterns constructed by the shogi crystal of p2mm symmetry and the bishop, silver general, gold general, rook or king, the outside and base neighborhood are included in their control.

Second we consider the shogi crystal of p2 symmetry.
As shown in Fig. \ref{Fig16}, the shogi patterns constructed by the shogi crystal of p2 symmetry and any kinds of the shogi pieces except the knight and the pawn satisfy the nearly complete neighborhood control condition.
Also, the shogi patterns constructed by the shogi crystal of p2 symmetry and the knight or pawn do not satisfy the nearly complete neighborhood control condition.
One can show in the shogi patterns constructed by the shogi crystal of p2 symmetry and the knight or pawn, a part of the spaces in the outside neighborhood is not included in their control.
Also, in the shogi patterns constructed by the shogi crystal of p2 symmetry and the lance, bishop, or silver general, the outside (base) neighborhood is (is not) included in their control.
In addition, in the shogi patterns constructed by the shogi crystal of p2 symmetry and the the gold general, rook or king, the outside and base neighborhood are included in their control.

Third we consider the shogi crystal of p1m1 symmetry.
As shown in Fig. \ref{Fig17}, the shogi patterns constructed by the shogi crystal of p1m1 symmetry and any kinds of the shogi pieces except the knight and the pawn and the lance satisfy the nearly complete neighborhood control condition.
Also, the shogi patterns constructed by the shogi crystal of p1m1 symmetry and the knight, pawn, or lance do not satisfy the nearly complete neighborhood control condition.
One can show in the shogi pattern constructed by the shogi crystal of p1m1 symmetry and the knight, pawn, or lance, a part of the spaces in the outside neighborhood is not included in their control.
In addition, in the shogi patterns constructed by the bishop, silver general, gold general, rook, or the king, the outside (base) neighborhood is (is not) included in their control.

Fourth we consider the shogi crystal of p11m symmetry.
As shown in Fig. \ref{Fig18}, the shogi patterns constructed by the shogi crystal of p11m symmetry and the silver general, gold general, rook, or king satisfy the nearly complete neighborhood control condition.
In addition, the shogi patterns constructed by the shogi crystal of p11m symmetry and the knight, pawn, lance, or bishop do not satisfy the nearly complete neighborhood control condition.
Note that in the shogi patterns constructed by the shogi crystal of p11m symmetry and the knight, pawn, lance, or bishop, a part of the spaces in the outside neighborhood is not included in their control.
In addition, in the shogi pattern constructed by the shogi crystal of p11m symmetry and the silver general, the outside and base (inner) neighborhood are (is not) included in the control.
Also, in the shogi pattern constructed by the shogi crystal of p11m symmetry and the gold general, the outside and inner (base) neighborhood are (is not) included in the control.
In addition, in the shogi pattern constructed by the shogi crystal of p11m symmetry and the rook or king, the outside, inner, and base neighborhood are included in the control.

Fifth we consider the shogi crystal of p2mg symmetry.
As shown in Fig. \ref{Fig19}, the shogi patterns constructed by the shogi crystal of p2mg symmetry and the gold general, rook, or king satisfy the nearly complete neighborhood control condition.
In addition, the shogi patterns constructed by the shogi crystal of p2mg symmetry and the knight, pawn, lance, bishop, or silver general do not satisfy the nearly complete neighborhood control condition.
One can show in the shogi patterns constructed by the shogi crystal of p2mg symmetry and the knight, pawn, lance, bishop, or silver general, a part of the outside neighborhood is not included in the control.
Also, in the shogi pattern constructed by the shogi crystal of p2mg symmetry and the gold general or rook, the outside (base) neighborhood is (is not) included in the control.
Also, in the shogi pattern constructed by the shogi crystal of p2mg symmetry and the king, the outside and base neighborhood are included in the control.

Sixth we consider the shogi crystal of p1 symmetry.
As shown in Fig. \ref{Fig20}, the shogi patterns constructed by the shogi crystal of p1 symmetry and the rook or king satisfy the nearly complete neighborhood control condition.
In addition, the shogi patterns constructed by the shogi crystal of p1 symmetry and the knight, pawn, lance, bishop, silver general, or gold general do not satisfy the nearly complete neighborhood control condition. 
One can show in the shogi patterns constructed by the shogi crystal of p1 symmetry and the knight, pawn, lance, bishop, silver general, or gold general, a part of the outside neighborhood is not included in the control.
Also, in the shogi pattern constructed by the shogi crystal of p1 symmetry and the rook or king, the outside (base) neighborhood is (is not) included in the control.

Seventh we consider the shogi crystal of p11g symmetry.
As shown in Fig. \ref{Fig21}, the shogi pattern constructed by the shogi crystal of p11g symmetry and the king satisfies the nearly complete neighborhood control condition.
Also, the shogi patterns constructed by the shogi crystal of p11g symmetry and the other kinds of shogi pieces do not satisfy the nearly complete neighborhood control condition.
Note that in the shogi patterns constructed by the shogi crystal of p11g symmetry and the shogi crystal of p11g symmetry and the knight, pawn, lance, bishop, silver general, gold general, or rook, a part of the outside neighborhood is not included in the control.
Also, in the shogi pattern constructed by the shogi crystal of p11g symmetry and the king, the outside neighborhood is included in the control.


If the movements of the shogi pieces were slightly different, Fig. \ref{Fig14} would not show a clear correspondence.
For example, if the lance can move to the back, or the silver general can move to the left or right, or the knight can move the same way as the knight of chess, the correspondence is broken.
Thus, shogi crystals are clearly associated with the movement of the shogi pieces to some extent.
Note that if the lance was a piece that could advance only two squares forward, instead of as far as it could go forward, the same argument as in this paper would still hold, and that level of arbitrariness would remain.
However, given the complexity of shogi piece movements, this correspondence may well be worthwhile.

\newpage
\section{Conclusion}
In this paper, we studied the mathematical aspects of the movements of shogi pieces moving away from the study of winning and losing in shogi. 
First we introduced the concept of the shogi pattern which is a geometric pattern constructed by shogi pieces and proposed the condition concerning the neighborhood for shogi patterns, nearly complete neighborhood control condition.
Next we proposed the shogi crystals that are based on the shape of the shogi pieces and whose symmetries are characterized by the Frieze groups.
Then we obtained the link between the shogi crystals and the movement of the shogi pieces through the condition of the nearly complete neighborhood control condition.
This paper may be developed in many directions. 
For example, as mentioned in the introduction, ancient shogi games included more types of pieces, and it would be interesting to address the mathematical structure of those as well.
Also, while this paper mainly dealt with shogi patterns consisting of one kind of shogi pieces, it would be fascinating to study shogi patterns consisting of multiple kinds of shogi pieces in more detail.
Furthermore, it would be possible to investigate the mathematical structure of the shogi pieces using not only the Frieze groups but also various geometric patterns such as the wallpaper group and space group.


\end{document}